\documentclass[12pt,a4paper]{article}

\usepackage{a4wide,amssymb,amsmath,amsthm,xspace,epsfig, amsfonts}

\begin{document}

\newcommand{\N}{\mathbb{N}}
\newcommand{\R}{\mathbb{R}}
\newcommand{\Z}{\mathbb{Z}}
\newcommand{\Q}{\mathbb{Q}}
\newcommand{\C}{\mathbb{C}}
\newcommand{\PP}{\mathbb{P}}

\newcommand{\LL}{\Bbb L}
\newcommand{\OO}{\mathcal{O}}

\newcommand{\esp}{\vskip .3cm \noindent}
\mathchardef\flat="115B

\newcommand{\lev}{\text{\rm Lev}}

\def\ut#1{$\underline{\text{#1}}$}
\def\CC#1{${\cal C}^{#1}$}
\def\h#1{\hat #1}
\def\t#1{\tilde #1}
\def\wt#1{\widetilde{#1}}
\def\wh#1{\widehat{#1}}
\def\wb#1{\overline{#1}}

\def\restrict#1{\bigr|_{#1}}

\def\hu#1#2{\mathsf{U}_{fin}\bigl({#1},{#2}\bigr)}
\def\ch#1#2{\left(\begin{array}{c}#1 \\ #2 \end{array}\right)}

\newtheorem{lemma}{Lemma}[section]

\newtheorem{thm}[lemma]{Theorem}

\newtheorem{defi}[lemma]{Definition}
\newtheorem{conj}[lemma]{Conjecture}
\newtheorem{cor}[lemma]{Corollary}
\newtheorem{prop}[lemma]{Proposition}
\newtheorem*{prob}{Problem}
\newtheorem{q}[lemma]{Question}
\newtheorem*{rem}{Remark}
\newtheorem{examples}[lemma]{Examples}
\newtheorem{example}[lemma]{Example}

\title{Notes on linearly H-closed spaces and od-selection principles}
\date{\empty}
\maketitle

\abstract{A space is called linearly H-closed if and only if any chain cover possesses a dense member. 
This property lies strictly between feeble compactness and H-closedness.
While regular H-closed spaces are compact, there are non-compact 
linearly H-closed spaces which are even collectionwise normal and Fr\'echet-Urysohn.
We give examples in other classes, and ask whether there is a first countable normal linearly H-closed non-compact space in {\bf ZFC}.
We show that {\bf PFA} implies a negative answer if the space is moreover either locally separable or 
both locally compact and locally ccc. An
Ostaszewski space (built with $\diamondsuit$) is an example which is even perfectly normal.
We also investigate Menger-like properties for the class of od-covers, that is, covers whose members are open and dense.}

\section{Introduction}

This note is mainly about a property (to our knowledge not 
investigated before) we decided to call linear H-closedness, which lies strictly between H-closedness and feeble compactness.
Since it came up while investigating simple instances of od-selection properties (see below), 
and all have a common `density of open sets' flavor,
we included a section about this latter topic although they are not related more than on a superficial level.

By `space' we mean `topological space'.  
We take the convention that `regular' and 'normal' imply `Hausdorff'.
A {\em cover} of a space always means a cover by open sets, and 
a cover is a {\em chain cover} if it is linearly ordered by the inclusion relation.
In any Hausdorff space (of cardinality at least $2$), each point has a non-dense neighborhood,
and thus the space has the property of possessing a cover by open non-dense sets.
But the chain-generalization of this property may fail.

\begin{defi}
  A space $X$ is linearly H-closed if and only if any chain cover has a member which is dense in $X$ (or equivalently if and only if any chain cover has a finite subfamily with a dense union).
\end{defi}  

Recall that a Hausdorff space 
any of whose covers has a finite subfamily with a dense union is called {\em H-closed}, whence the name `linearly H-closed'.
While H-closed regular spaces are compact 
(see \cite[Corollary 4.8(c)]{PorterWoodsBook} for a simple proof),
there are plenty of Tychonoff linearly H-closed non-compact spaces, 
perhaps the most simple being the
Tychonoff plank
(see Example \ref{lemmalinHbigL}).
We will give examples in various classes such as first countable, normal, collectionwise normal, etc,
but while there are consistent examples of non-compact perfectly normal first countable linearly H-closed spaces, we 
were unable to determine whether a non-compact first countable normal linearly H-closed space exists in {\bf ZFC} alone.
A partial result is that {\bf PFA} prevents such a space from existing if it is moreover either locally separable or 
both locally compact and
locally ccc (see Theorem \ref{prophersep}).
These results are contained in Section \ref{sec1}.

In Section \ref{sec2}, we investigate
Menger-like properties for {\em od-covers} of topological spaces, that is, covers whose members are open and dense. 
In our short study, we show in particular that the class of non-compact spaces satisfying $\hu{\OO}{\Delta}$ does contain some 
Hausdorff spaces
but no regular space, and that a separable space satisfies $\hu{\Delta}{\OO}$ if and only if it satisfies $\hu{\OO}{\OO}$, where $\Delta$
is the class of od-covers. We defer the definitions of $\hu{\mathcal{A}}{\mathcal{B}}$
until Section \ref{sec2}. 
Research on selection principles (such as Menger-like properties)
currently flourishes and sees an impressive flow of new results
(see for instance \cite{Scheepers:survey, Tsaban:survey} for surveys about recent activity in the field).
Since the author is not an expert on the subject and admits to feeling
a bit lost in its numerous subtleties, we shall
content ourselves with a humble introduction to the class of od-covers and derive only basic properties. 

For convenience, we now give a grouped definition:
the {\em (od-)[linear-]Lindel\"of number} ($odL(X)$) [$\ell L(X)$] $L(X)$ of a space $X$ is the smallest cardinal $\kappa$ such
that any (od-)[chain] cover of $X$ has a subcover of cardinality $\le\kappa$.
A space is {\em od-compact} if and only if any od-cover has a finite subcover, and we
define similarly {\em od-Lindel\"of} , {\em linearly-Lindel\"of}, etc.
We do not assume separation axioms in any of these properties. It happens that
the od-Lindel\"of number and the Lindel\"of number almost always coincide, the only exception is when
the space contains a `big' clopen discrete subspace.
See Section \ref{sec2} (especially Theorem \ref{ehouais}) for details and remarks about
the ignorance of past results.
For the information of the reader, we note that our definitions above of 
$L(X),odL(X)$ and $\ell L(X)$ are different from the ones we found convenient to give in \cite{mesziguesod},
where, for example, $L(\R) = \omega_1$, not $\omega$.

{\em Acknowledgments.} The idea of looking at selection principles for od-covers has been given to us by B. Tsaban.
The desire to complete this paper came from the need of writing up a Curriculum Vit\ae \
(for which this paper turned out to be irrelevant) M. Lazeyras asked for.
We thank both of them, as well as Z. and C. Petrini (for personal reasons they know about).
We also thank the anonymous referee of this paper for his/her numerous suggestions and corrections. Some of our 
results are actually due to him/her (and quoted such).


\section{Linearly H-closed spaces}\label{sec1}

In this section, each space is assumed to be Hausdorff, even though that property 
is not needed for every assertion, and we will repeat the assumption often (for clarity).
Any chain cover possesses a subcover indexed by a regular cardinal and for simplicity we will always use such indexing.
It is immediate that the continuous image of a linearly H-closed space is linearly H-closed.
Our first lemma is almost trivial.

\begin{lemma}\label{lemma:trivialone}
  A space is linearly H-closed if and only if any infinite cover of it has a subfamily of strictly smaller cardinality with a dense union.
\end{lemma}
\proof 
  Given
  a chain cover indexed by a regular cardinal, a subcover of strictly smaller cardinality is contained in some member, so the latter implies the former.
  If $X$ is linearly H-closed,
  given a cover $\{U_\alpha\,:\,\alpha\in\kappa\}$, then the sets $V_\alpha=\cup_{\beta<\alpha}U_\beta$ form a chain cover and 
  some $V_\alpha$ is dense.
\endproof

It is well known that a space is H-closed if and only if any open filter base on $X$ 
(that is, a filter base containing only open subsets of $X$) has an adherent point.
See for instance (the proof of) \cite[Prop. 4.8(b) (2) $\Leftrightarrow$ (3)]{PorterWoodsBook}.
The referee pointed out to us that a similar result holds for linear H-closedness.
By a {\em chain filter base} we mean an open filter base which is linearly ordered by the inclusion relation.
The proof we just mentioned
can be easily adapted to show the following.

\begin{lemma}\label{lemma:ref1}
  A space $X$ is linearly H-closed if and only if any chain filter base on $X$ has an adherent point.
\end{lemma}

Likewise, the following result (also suggested by the referee) can be proved as Proposition 4.8(e) in \cite{PorterWoodsBook}.

\begin{lemma}\label{lemma:ref2}
  If $X$ is linearly H-closed and $U$ is open, then $\wb{U}$ is linearly H-closed.
\end{lemma}

\begin{figure}[h]
  \begin{center}
  \epsfig{figure=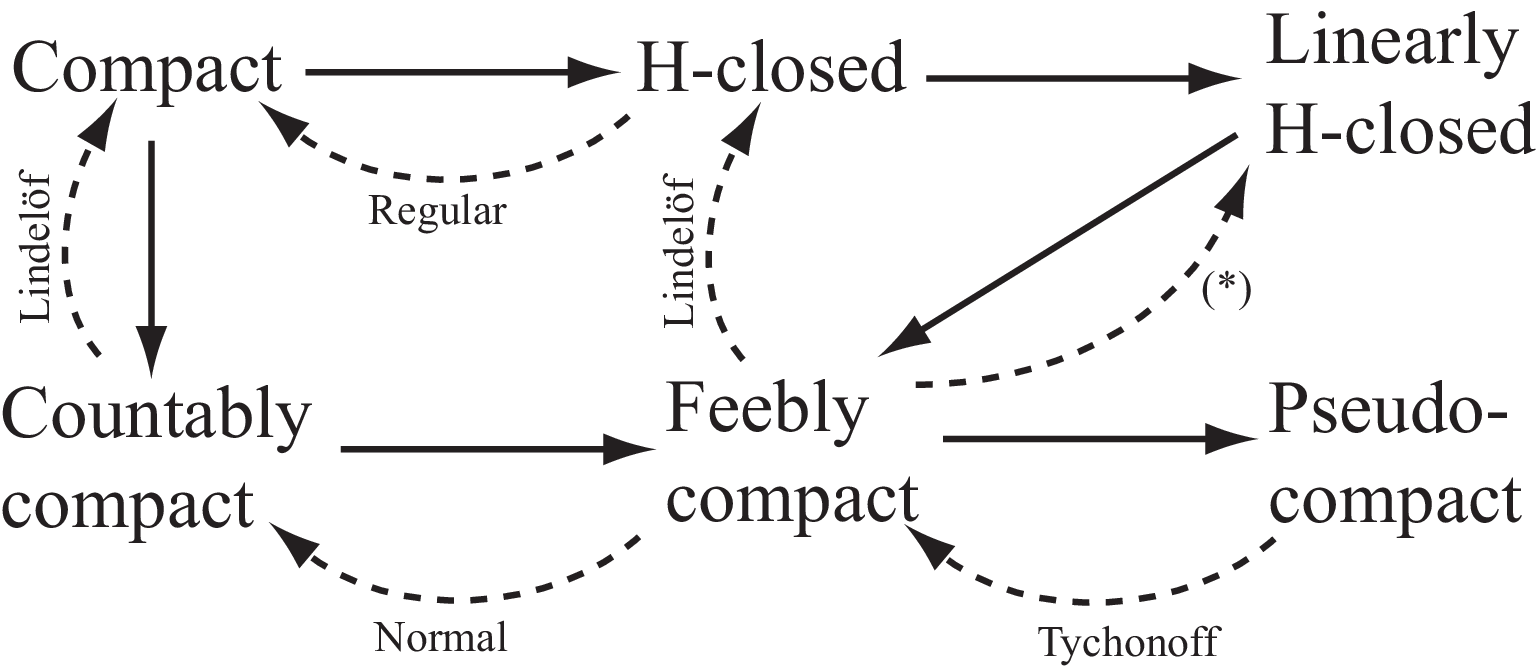, width=10cm}
  \caption{Some implications for Hausdorff spaces.}\label{fig:1}
  \end{center}
\end{figure}

However, not every closed subset of a linearly H-closed space is linearly H-closed;
see for instance Example \ref{lemmalinHbigL}.
Linear H-closedness is linked to other generalized
compactness properties, as seen in Figure 1 below. Plain straight arrows denote implications that 
hold for Hausdorff spaces (and most of them for any space) while additional properties (for instance those
written on their side) are needed for
those denoted by dotted curved arrows. 
Recall that a space is {\em feebly compact} if and only if every locally finite 
family of open sets is finite. This turns out
to be equivalent to ``every locally finite cover is finite'' and to 
``every countable cover of $X$ has a finite subfamily with a dense union'' (see \cite[Theorem 1.11(b)]{PorterWoodsBook}).
The term ``feebly compact'' is due to S. Marde\v sic and P. Papi\'c (see \cite[p. 902]{Stone:1960}).
A space is {\em pseudocompact} if and only if any continuous real valued function on it is bounded.
All implications in Figure \ref{fig:1} are classical except linearly H-closed $\longrightarrow$ feebly compact
and its converse whose proofs are given in Lemma \ref{lemmacccfc}. An example of Condition (*) is given in the statement of the lemma.

We decided to state this lemma in an almost absurd amount of generality, so we need some definitions.
The good news is that more readable corollaries do follow quite easily.
Given an infinite cardinal $\kappa$,
a space is {\em initially $\kappa$-[linearly] Lindel\"of} if and only if any open [chain] cover of cardinality $\le\kappa$ has a countable subcover. 
Notice that any space is initially $\omega$-[linearly] Lindel\"of.
The weak Lindel\"of number $wL(X)$ of a space $X$ is the least cardinal $\kappa$ such that
any open cover of $X$ has a subfamily of cardinality $\le\kappa$ whose union is dense.
Notice that if $Y\subset X$ is dense, then $wL(X)\le wL(Y)$ and if $Y$ is feebly compact then so is $X$.

\begin{lemma}\label{lemmacccfc} \ \\
    (1) A linearly H-closed space is feebly compact.\\
    (2) Let $X$ be a Hausdorff space, $Y\subset X$ be dense in $X$, and $\kappa$ be an infinite cardinal.
        Assume that 
        $wL(X)\le \kappa$ and that $Y$ is both initially $\kappa$-linearly Lindel\"of and feebly compact. 
    Then $X$ is linearly H-closed.
\end{lemma}
\proof
  \ \\
  (1) Given a countable cover $\mathcal{U}=\{U_n\,:\,n\in\omega\}$ of a linearly H-closed $X$, set $V_n=\cup_{m\le n}U_m$.
      Then $V_n$ is dense for some $n$, and the result follows.\\
  (2) 
   Let $\mathcal{U}=\{U_\alpha\,:\,\alpha\in\lambda\}$ be an infinite chain cover of $X$, with $\lambda$ a regular cardinal. 
   Assume first that $\lambda\le\kappa$.
   There is thus a countable subfamily that covers $Y$, and then some $U_\alpha$ is dense in it by feeble compactness. 
   It follows that $U_\alpha$ is dense in $X$ as well.
   Now, suppose that $\lambda>\kappa$. Since $wL(X)\le\kappa$
   there is some subfamily of cardinality $\le\kappa<\lambda$ whose union is dense in $X$, 
   and by regularity of $\lambda$ its union is contained in some $U_\alpha$.
\endproof

A case not covered by this lemma is the following easy fact.

\begin{lemma}\label{lemmafacile}
   Let $X$ be a Hausdorff space containing a dense feebly compact linearly Lindel\"of subspace $Y$. Then $X$ is linearly H-closed.
\end{lemma}
\proof
   Given a chain cover of $Y$, linear Lindel\"ofness gives a countable subcover and then feeble compactness gives a finite subfamily
   of the subcover which contains a dense member.
\endproof

For a cardinal $\kappa$, a space is {\em $\kappa$-cc} (or {\em ccc} if $\kappa=\omega$) 
if and only if any disjoint collection of open sets has cardinality at most $\kappa$. A space
with a dense subset of cardinality $\kappa$ is obviously $\kappa$-cc. 

\begin{cor}\label{corccc}
  If $X$ is Hausdorff and possesses a dense feebly compact ccc subspace $Y$, then $X$ is linearly 
  H-closed.
\end{cor}
\proof
It is well known that a $\kappa$-cc space has weak Lindel\"of number $\le\kappa$; hence $wL(Y)\le\omega$.
Invoking the vacuousness of the definition, $Y$ is also initially $\omega$-Lindel\"of, and the conditions of Lemma \ref{lemmacccfc} (2)
are thus fulfilled. 
\endproof

Recall that a space is {\em perfect} if and only if any closed subset is a $G_\delta$.

\begin{cor}\label{corperfect}
  Let $X$ be a feebly compact regular perfect space. Then $X$ is first countable and linearly H-closed.
\end{cor}
\proof
In paragraph (b) on page 378 of \cite{Glicksberg:1959}, I. Glicksberg (using different terminology) provided a proof that
a $G_\delta$ point in a regular feebly compact space has a countable neighborhood base.
For another proof, see
Lemma 2.2 in \cite{PorterWoods:1984}. Moreover, Lemma 2.3 in 
\cite{PorterWoods:1984} shows that
if each closed set in $X$ is a $G_\delta$, then
$X$ is ccc. 
\endproof

We can use Lemma \ref{lemmacccfc} to obtain simple examples:

\begin{example}\label{lemmalinHbigL}
   There are linearly H-closed Tychonoff spaces of arbitrarily high weak Lindel\"of number and cellularity.
\end{example}
\proof[Details]
   A very classical example: the Tychonoff plank of a regular cardinality.
   Let us recall the construction and its properties for convenience.
   Fix a regular cardinal $\kappa$. 
   Let $X$ be the subspace of the product $(\kappa^+ +1)\times(\kappa +1)$ obtained by removing the point $\{\langle\kappa^+,\kappa\rangle\}$.
   Each ordinal is given the order topology.
   \\
   As a subspace of a compact space, $X$ is Tychonoff. 
   The cellularity of $X$ is at least $\kappa^+$ since $\{\alpha\}\times(\kappa+1)$ for successor $\alpha\in\kappa^+$ 
   is a disjoint collection of open subsets. 
   The cover $\{\alpha\times(\kappa+1)\,:\,\alpha\in\kappa^+\}\cup\{(\kappa^++1)\times\beta\,:\,\beta\in \kappa\}$
   shows that $wL(X)\ge\kappa$.
   Since $(\kappa^+ +1)\times \kappa$ is the union of $\kappa$ compact sets and is dense in $X$, $wL(X)\le\kappa$.
   Recall that $\kappa^+$ with the order topology is initially $\kappa$-compact, and so is its product with the 
   compact space $(\kappa+1)$ (see, e.g., Theorem 2.2 in \cite{Stephenson:1984}). Thus the dense subset
   $Y=\kappa^+\times(\kappa+1)$ is in particular feebly compact and initially $\kappa$-Lindel\"of.
   This implies that $X$ is linearly H-closed by Lemma \ref{lemmacccfc} (2). 
\endproof

Of course, these spaces are not first countable. Let us give more elaborate
examples. All are `classical' spaces which happen to be linearly H-closed.
In the following, we refer to \cite{Vaughan:1990} for the definitions of the `small'
uncountable cardinals $\mathfrak{p},\mathfrak{b}$, but recall that $\omega_1\le\mathfrak{p}\le\mathfrak{b}\le 2^{\aleph_0}$ and
that each inequality may be strict. The diamond axiom $\diamondsuit$ implies the continuum hypothesis {\bf CH} and is defined in any book on set theory.

\begin{examples}\label{exBell}\ 
   There are linearly H-closed non-compact spaces with the following additional properties:\\
   (a) (Bell) First countable, Tychonoff, Lindel\"of number $\omega_1$.\\
   (b) (Isbell-Mr\'owka) First countable, locally compact (and thus Tychonoff), perfect.\\
   (c)  
         (Franklin and Rajagopalan, in effect) {\rm ($\mathfrak{p}=\omega_1$)} First countable, locally compact, normal.\\
   (d) (Ostaszewksi) {\rm ($\diamondsuit$)} First countable, locally compact, perfectly normal.\\
   (e) Frechet-Urysohn, collectionwise normal.\\
\end{examples}

\proof[Details]
Linear H-closedness follows from Corollary \ref{corccc} in each case except (b) where Corollary \ref{corperfect} is used.
\\
(a) M.G. Bell \cite[Example 1]{Bell:1989} constructed a first countable countably compact ccc (non-separable) Tychonoff space $X$.
Since $X$ is an increasing union of $\aleph_1$-many compact spaces, it has Lindel\"of number $\omega_1$. 
\\
(b) The space $\Psi$, due independently to J. Isbell's and S. Mr\'owka (see e.g. \cite[Exercise 5I]{GillmanJerisonBook} or \cite{Mrowka:1954}), 
is first countable, perfect,
Tychonoff and feebly compact. This space is not countably compact, and thus
non normal.\\
(c) 
Franklin and Rajagopalan introduced a class of spaces called $\gamma\N$ spaces, which consist of a dense 
discrete countable set to which is `attached' a copy of $\omega_1$
in such a way that the space is locally compact and normal, with various additional properties depending on how the attachment is done.
The constructions were later simplified and generalized by van Douwen, Nyikos and Vaughan, and a version of 
$\gamma\N$ which is countably compact
and first countable can be built if and only if $\mathfrak{p}=\omega_1$ (see for instance \cite{Nyikos:countablycompact1986}, Theorem 2.1 and Example 3.4,
or \cite{Nyikos:HereditaryNormality}).
\\
(d) The celebrated Ostaszewski's space \cite{Ostaszewski:1976}: a first countable,
perfectly normal, hereditarily separable,
countably compact, locally compact, non-compact space built with $\diamondsuit$.\\
(e) The sigma-product of $2^{\omega_1}$, i.e. the subspace of the compact space $2^{\omega_1}$ where at most countably
many coordinates have value $1$, is collectionwise normal, Frechet-Urysohn, countably compact and ccc
(see for instance H. Brandsma's answer on the MathOverflow question \cite{MO79021}).
\endproof

More than {\bf ZFC} is necessary for the construction in (c); see Theorem \ref{thmhersep} below.
Bell's space in (a) cannot be shown to be locally compact in {\bf ZFC} by Theorem \ref{thmcccsep}.
It is also not separable, and no separable regular example with Lindel\"of number $\omega_1$ can be found in {\bf ZFC}, 
as the next lemma shows.

\begin{lemma}\label{lemmafcb}
  A first countable separable linearly H-closed Hausdorff space of Lindel\"of number $<\mathfrak{p}$ is H-closed 
  (and thus compact if regular).
\end{lemma}
Notice the similarity with the fact (proved in \cite{Hechler:1975}) that
a regular separable countably compact space of Lindel\"of number $<\mathfrak{p}$ is compact. 
\proof
A first countable separable space has
countable $\pi$-weight, as easily seen. Since $X$ is linearly H-closed, it is feebly compact.
A feebly compact space with countable $\pi$-weight and Lindel\"of number $<\mathfrak{p}$ is H-closed
(Lemma 3.1 in \cite{PorterWoods:1984}).
\endproof

Likewise, Example \ref{exBell} (d) cannot be constructed in {\bf ZFC + CH} alone. 

\begin{lemma}
  It consistent with {\bf ZFC} (and even with {\bf ZFC + CH}) that a perfectly normal linearly H-closed
  space is compact. In particular, it follows from {\bf MA + $\neg$CH}.
\end{lemma}
\proof
  A linearly H-closed normal space is countably compact, Weiss \cite{Weiss:1978} showed that {\bf MA + $\neg$CH} 
  implies that a countably compact regular perfect space is compact, and Eisworth \cite{Eisworth:2001} showed that this
  latter
  result is compatible with {\bf CH}.
\endproof

\begin{q}
   Is there a normal first countable linearly H-closed non-compact space in {\bf ZFC}?
\end{q}

The following theorem is a partial answer.

\begin{thm}\label{prophersep}
   {\rm ({\bf PFA})}
   Let $X$ be a normal linearly H-closed space. 
   If either
   (a) $X$ is countably tight and locally separable, or 
   (b) $X$ is first countable, locally compact and locally ccc, then $X$ is compact.
\end{thm}

Our use of {\bf PFA} is indirect. Indeed, we only need two of its classical consequences.
Recall that {\bf PFA} implies {\bf MA + $\neg$CH}.

\begin{thm}\label{thmhersep}
   (Balogh-Dow-Fremlin-Nyikos, \cite{Balogh-Dow_Fremlin-Nyikos:1988}, Corollary 2)\\
    {\rm ({\bf PFA})}
   Every separable, normal, countably tight, countably compact space is compact.
\end{thm}

\begin{thm}\label{thmcccsep}
   (Hajnal-Juh\'asz \cite{HajnalJuhasz:1971}) \\
   {\rm({\bf MA + $\neg$CH})}
   Every first countable, locally compact, ccc space is separable.
\end{thm}

\proof[Proof of Theorem \ref{prophersep}]
  \ \\
  (a) If $X$ is not compact, we will show that
      it is possible to
      define open subsets $U_\alpha\subset X$, for each $\alpha<\omega_1$, such that $\wb{U_\beta}\subsetneqq U_\alpha$ whenever 
      $\beta <\alpha$. Then $Y=\cup_{\alpha<\omega_1}U_\alpha$ is a clopen subset of $X$. Openness is immediate. To see that it is closed,
      notice that given a point $x\in\wb{Y}$ by countable tightness 
      there is a countable subset of $Y$ having $x$ in its closure. But a countable subset of $Y$
      is contained in some $U_\alpha$, so $x\in\wb{U_\alpha}\subset U_{\alpha+1}\subset Y$.     
      It follows that $X$ is not linearly H-closed, since no member of the chain cover 
      $\{(X-Y)\cup U_\alpha\,:\,\alpha\in\omega_1\}$ is dense in $X$. 
      \\
      To find $U_\alpha$, we proceed by induction. Each will be a separable open subset of $X$.
      Let $U_0$ be any such open separable subset.
      Assume that $U_\beta$ is defined for each $\beta<\alpha$.
      Recall that by normality and linear H-closedness $X$ is countably compact.
      Thus, $Z=\wb{\cup_{\beta<\alpha}U_\beta}$, being separable, is compact by
      Theorem \ref{thmhersep}. If $Z=X$, then $X$ is compact.
      Otherwise choose a point $x\not\in Z$, cover $\{x\}\cup Z$
      by open separable sets and take the union of a finite subcover to obtain
      a separable $U_{\alpha}$ properly containing $Z$. In particular $\wb{U_{\beta}}\subset U_{\alpha}$ for all $\beta <\alpha$. 
      This defines $U_\alpha$ for each $\alpha<\omega_1$ with the required properties.\\
  (b) We proceed as in (a), defining $U_\alpha$ to be ccc with compact closure.
      The successor stages are the same. If $\alpha$ is limit then $\wb{\cup_{\beta<\alpha}U_\beta}$, having a dense ccc subspace, is ccc.
      By Theorem \ref{thmcccsep}, it is separable under {\bf MA + $\neg$CH} and thus compact under {\bf PFA}.
\endproof

Note: Theorem 5.4 in \cite{Nyikos:countablycompact1986} seems to indicate that there are models of {\bf MA + $\neg$CH} or even {\bf PFA$^-$} 
with separable, locally compact, locally countable, countably compact, countably tight normal spaces,
but we do not know to which spaces this assertion refers.
The referee kindly 
informed us that the preprint \cite{NyikosForcingCompact} by P. Nyikos, where theses spaces were probably described, was never published.

We now briefly investigate how 
far a first countable linearly H-closed space is from being sequentially compact and show in
Lemma \ref{propno} below that there are restrictions on the Lindel\"of number. 
(The result seems well known, see the remarks before
Problem 359 in \cite{Vaughan:1990}, but we include the proof for completeness.)
We first need some vocabulary. A collection of subsets of $X$ is a {\em discrete collection} if each point of $X$ 
possesses a neighborhood intersecting at most one member of the collection. This implies that given any subcollection,
the union of the closures of its members is closed.
A space satisfies the condition {\em wD} if given any infinite closed discrete subspace $D$ of $X$, there
is an infinite $D'\subset D$ which expands to a discrete collection of open sets, that is,
for each $x\in D'$ there is an open $U_x\ni x$ such that $\{U_x\,:\,x\in D'\}$ is a discrete collection.

\begin{lemma}\label{propno}
  A regular, first countable, feebly compact space is either countably compact or has Lindel\"of number $\ge\mathfrak{b}$.
\end{lemma}

\proof  
  Let $X$ be a regular first countable space whose Lindel\"of number is $<\mathfrak{b}$ and suppose that it is not countably compact.
  Let thus $\{x_n\in X\,:\,n\in\omega\}$ be an infinite closed discrete subset.
  A regular first countable space with Lindel\"of number $<\mathfrak{b}$
  satisfies {\em wD} (see 3.6 \& 3.7 in \cite{Nyikos:1990}).
  Let thus $E\subset\omega$ be infinite and $U_n\ni x_n$ ($n\in E$) be open such that $\{\wb{U_n}\,:\,n\in E\}$ is discrete.
  In particular $\{U_n\,:\,n\in E\}$ is an infinite locally finite family of open sets, which is impossible in a feebly compact space.
\endproof

We close this section with two results due to the referee who kindly gave us permission to include them in this note.
Firstly, notice that by continuity of the projections, if a product of spaces is linearly H-closed then each factor space is linearly H-closed.
But the converse may fail:

\begin{prop}\label{prop:ref3}
   There is a linearly H-closed space $G$ such that $G\times G$ is not linearly H-closed.
\end{prop}
\proof
It is well known (e.g., see \cite[Ex. 9.15]{GillmanJerisonBook}) that there exists a subspace $G$ of the Stone-\v Cech compactification
of the integers $\beta\omega$ such that $\omega\subset G$ (hence $G$ is separable), and $G$, but not $G\times G$, is feebly compact.
Thus, $G$ is linearly H-closed by Corollary \ref{corccc} while $G\times G$ is not.
\endproof

However, the following holds:

\begin{prop}\label{prop:ref4}
   If $X$ is H-closed and $Y$ is linearly H-closed, then $X\times Y$ is linearly H-closed.
\end{prop}
\proof
We use the characterization of linear H-closedness given by Lemma \ref{lemma:ref1}.
Let $\mathcal{U}$ be a chain filter base on $X\times Y$. Since the projection on the $Y$ factor $\pi_Y$ is open, 
$\{\pi_Y(U)\,:\,U\in\mathcal{U}\}$ is a chain filter base on $Y$ and hence has an adherent point $y\in Y$.
Let 
$$
  \mathcal{P} = \{(X\times W)\cap U\,:\, U\in\mathcal{U}\text{ and }W\text{ is an open neighborhood of } y\}.
$$
Then $\mathcal{P}$ is an open filter base, and hence 
$\{\pi_X(V)\,:\,V\in\mathcal{P}\}$ must have an adherent point $x\in X$.
For every neighborhood $V\subset X$ and $W\subset Y$ of $x,y$, respectively, and every $U\in\mathcal{U}$, we have 
by construction $U\cap(V\times W)\not=\varnothing$.
It follows that $\langle x,y\rangle$ is an adherent point of $\mathcal{U}$.
\endproof


\section{Od-selection properties}\label{sec2}

No separation axiom is assumed in this section.
Allow us first a remark about the od-Lindel\"of number.
The author proved in 
\cite{mesziguesod} that a $T_1$ space is od-compact if and only if the subspace of non-isolated points is compact,
and that a $T_1$ space with od-Lindel\"of number $\le\kappa$ either has a closed discrete subset of cardinality $>\kappa$,
or $\ell L(X)\le\kappa$ holds whenever $\kappa$ is regular.
We made the remark that since the methods were elementary, it would not
be a surprise if similar results we were unaware of had appeared elsewhere.
It was indeed the case: Mills and Wattel \cite{MillsWattel} had 
shown that a $T_1$ space without isolated points with $odL(X)\le\kappa$ satisfies $L(X)\le\kappa$ as well, which
is much stronger (the compact case is actually due to Kat\v etov in 1947 \cite{Katetov:1947}).
Blair \cite{Blair:1983} later improved their proof. (Both papers
actually deal with $[\kappa,\lambda]$-compactness.) 
We show below that a very small modification of Blair's proof yields the following.

\begin{thm}(Mills--Wattel and Blair)\label{ehouais} 
  Let $\kappa$ be an infinite cardinal. 
  Let $X$ be a $T_1$ space with $odL(X)\le\kappa$. Then either $X$ contains a clopen discrete subset
  of cardinality $>\kappa$, or $L(X)\le\kappa$.
  Moreover, the subspace of non-isolated points of $X$ has Lindel\"of number $\le\kappa$.
\end{thm}
\proof[Quick proof, following Blair]
It is easy to see that a space has od-Lindel\"of number $\le\kappa$ if and only if any closed nowhere dense subset
has Lindel\"of number $\le\kappa$.
Let $\mathcal{U}$ be an open cover of $X$.
Let $\mathcal{W}$ be a maximal family of disjoint open sets such that each member of $\mathcal{W}$ is contained
in a member of $\mathcal{U}$. Then $\cup\mathcal{W}$ is dense. We may thus cover $X-\cup\mathcal{W}$ by 
a subfamily $\mathcal{V}\subset\mathcal{U}$ of cardinality $\le\kappa$. Take one point in each member
of $\mathcal{W}$ which is not
entirely covered by $\cup\mathcal{V}$.
This defines a closed discrete subset of $D\subset X$. 
Let $D_0=\{d\in D\,:\,d\text{ is isolated in }X\}$. 
Then $D-D_0$ is nowhere dense and hence of cardinality at most $\kappa$. 
It follows that at most $\kappa$ members of $\mathcal{U}$ cover $\cup\{W_d\,:\,d\in D-D_0\}$, where
$W_d$ is the unique member of $\mathcal{W}$ containing $d\in D$.
The uncovered part of $X$ is now contained in $\cup\{W_d\,:\,d\in D_0\}$.
Then either $|D_0|\le\kappa$, in which case we add $\le\kappa$ members of $\mathcal{U}$ to complete the subcover,
or
$|D_0|>\kappa$ and $X$ contains a clopen discrete subset
of cardinality $>\kappa$.
The `moreover' part follows easily from, e.g., Lemma 4.8 in \cite{mesziguesod}.
\endproof

For other results in the same spirit, see \cite{GartsideGreenwoodMcIntire}.
Let us now turn to selections properties. In what follows,
$\OO,\Delta$ respectively mean the collection of covers and od-covers of some topological space
which will be clear from the context. Recall that a cover is an {\em od-cover} if and only if every member is dense.
Given collections $\mathcal{A},\mathcal{B}$ of covers of a space $X$, we define the following property:
\begin{itemize}\item[\ ]
$\hu{\mathcal{A}}{\mathcal{B}}$~: 
For each sequence $\langle \mathcal{U}_n \colon n\in\omega\rangle$ of members of $\mathcal{A}$ which do not
have a finite subcover, there are finite $\mathcal{F}_n\subset \mathcal{U}_n$ such that $\{\cup \mathcal{F}_n\,:\,n\in\omega\} \in \mathcal{B}$.
\end{itemize}
Recall that the classical Menger property is (equivalent to) $\hu{\OO}{\OO}$, and that
$$
\sigma\text{-compact} \longrightarrow \hu{\OO}{\OO}\longrightarrow\text{ Lindel\"of}.
$$


\subsection*{The Property $\hu{\OO}{\Delta}$}

Let us first show the following simple lemma.

\begin{lemma}\label{cbhu}
        The following equivalences hold for any space $X$.
         \\
        (a) Lindel\"of $\&$ linearly H-closed $\longleftrightarrow$ Lindel\"of $\&$ H-closed, 
           \\
        (b) $
            \hu{\OO}{\OO} \, \& \text{ linearly H-closed } \longleftrightarrow 
           \hu{\OO}{\OO} \, \&  \text{ H-closed}  \longleftrightarrow  \hu{\OO}{\Delta}.
              $ 
    \\
    Moreover, the properties in (b) imply those in (a).
\end{lemma}
\proof 
  The moreover part is immediate since $\hu{\OO}{\OO}\longrightarrow$ Lindel\"of.\\
  (a) Immediate by Lemma \ref{lemma:trivialone}.\\
  (b) The leftmost equivalence follows from (a) by Lindel\"ofness. Let us prove the rightmost equivalence.
  For the direct implication, 
  let $\langle\mathcal{U}_n\rangle$ be a sequence of covers, and let $\mathcal{F}_n\subset\mathcal{U}_n$ be finite
  such that $\{\cup\mathcal{F}_n\,:\,n\in\omega\}$ is a cover of $X$. 
  By H-closedness, we can choose finite $\mathcal{G}_n\subset\mathcal{U}_n$
  such that $\cup\mathcal{G}_n$ is dense. Taking $\mathcal{F}_n\cup\mathcal{G}_n$ yields the result.
  For the converse implication,  $\hu{\OO}{\OO}$ trivially holds. We prove that $X$ is linearly H-closed and use the leftmost equivalence
  to obtain H-closedness. 
  Suppose that there is a chain cover $\mathcal{U}=\{U_n\,:\,n\in\omega\}$ without any dense member. 
  A finite union of members of $\mathcal{U}$, being contained in a member of $\mathcal{U}$, is therefore not dense,
  taking $\mathcal{U}_n =\mathcal{U}$ for all $n\in\omega$ gives a sequence of open covers violating $\hu{\OO}{\Delta}$. 
\endproof

The situation is then very simple for regular spaces:

\begin{prop} \label{propregular}
  The following properties are equivalent for regular spaces.
   \\
  (a) Lindel\"of $\&$ linearly H-closed,\\
  (b) $\hu{\OO}{\Delta}$,\\
  (c) Compact.
\end{prop}

\proof
  (b) $\rightarrow$ (a) 
  by Lemma \ref{cbhu} and (c) $\rightarrow$ (b) is trivial.
  Since
  a regular H-closed space is compact, (a) $\rightarrow$ (c) follows again by Lemma \ref{cbhu}. 
\endproof

We will show that both (a) $\rightarrow$ (b) and (b) $\rightarrow$ (c) may fail for Hausdorff spaces, that is,
we shall exhibit Hausdorff examples of Lindel\"of (linearly) H-closed spaces which do not satisfy $\hu{\OO}{\OO}$, 
and non-compact spaces satisfying $\hu{\OO}{\Delta}$.
Recall that a space $Y$ is an {\em extension} of the space $X$ if and only if 
$Y$ contains a copy of $X$ which is dense in $Y$.
H-closed extensions of Hausdorff spaces are well studied, see for instance
\cite{PorterWoodsBook}.  
The examples we describe below are very similar to the ones given in Chapter 7 of this book.
They can be seen as modifications of the half disk topology.

Let $X$ be a space equipped with two topologies $\tau,\rho$. Denote by $\wh{X}(\tau,\rho)$ the space 
whose underlying set is $X\times[0,1]$ topologized 
as follows. 
The topology on $X\times(0,1]$ is the product topology of $\tau$ and the usual metric topology on $(0,1]$.
Neighborhoods of $\langle x,0\rangle$ are then defined to be 
$U\times\{0\}\sqcup V\times (0,a)$ for $U\in\rho$, $V\in\tau$ with $x\in U\cap V$, and $0<a\le 1$.

\begin{lemma}\label{lemmawhX}
  Assume $\tau\subset\rho$, that is, $\rho$ is finer than $\tau$.\\
  (1) If $X$ is Hausdorff for $\tau$ (and thus for $\rho$), then so is $\wh{X}(\tau,\rho)$.\\
  (2) If $X$ is H-closed for $\tau$, then $\wh{X}(\tau,\rho)$ is H-closed.\\
  (3) $X$ is Lindel\"of for $\rho$ if and only if $\wh{X}(\tau,\rho)$ is Lindel\"of.\\
  (4) If $X$ is first countable for both $\tau$ and $\rho$, then so is $\wh{X}(\tau,\rho)$.
\end{lemma}

\proof
  Denote by $\tau\times\mu$ the product topology of $\tau$ on $X$ and the usual metric topology $\mu$ on $[0,1]$.
  Notice that the topology
  on $\wh{X}(\tau,\rho)$ is finer than $\tau\times\mu$ since $\tau\subset\rho$.
  \\
  (1) Immediate since $\tau\times\mu$ is Hausdorff.\\
  (2)
  A direct proof is not difficult, but let us give a more general argument suggested by the referee.
  Since the property H-closed is known to be productive (see e.g. \cite[Prop. 4.8(l)]{PorterWoodsBook}), 
  it follows that $Z = X\times[0,1]$, with topology $\tau\times\mu$,
  is an H-closed extension space of $Y = X\times (0,1]$. 
  $\wh{X}(\tau,\rho)$ is also an extension of $Y$ with the same underlying set as $Z$ and a finer topology.
  Moreover, $\wh{X}(\tau,\rho)$ and $Z$ have the same neighborhood filter trace on $Y$ in the sense
  of \cite[Def. 7.1(a)]{PorterWoodsBook}. Then Propositions 7.1(h) and 7.1(i) of \cite{PorterWoodsBook}
  imply that $\wh{X}(\tau,\rho)$ is H-closed.\\
  (3) The necessity is obvious since $X$ with the topology $\rho$ is a closed subspace of $\wh{X}(\tau,\rho)$.
      For the sufficiency, assume that $X$ is Lindel\"of for $\rho$. Then $X$ is Lindel\"of for $\tau$ as well.
      Since $X\times(0,1]$ with topology $\tau\times\mu$ is the product of a Lindel\"of space and a $\sigma$-compact space, it is Lindel\"of.
      It follows that $\wh{X}(\tau,\rho)=X\times\{0\}\cup X\times(0,1]$ is Lindel\"of.\\
  (4) Straightforward: a neighborhood basis for $\langle x,0\rangle$ is given by 
      $\{ U_n\times\{0\}\sqcup V_m\times (0,1/\ell)\,: \ell,m,n\in\omega,\,\ell>0\}$,
      where $U_n$, $V_n$ are local bases for $x$ in the $\rho$ and $\tau$ topologies.
\endproof

Notice that in most cases $\wh{X}(\tau,\rho)$ is not regular.

\begin{prop}\label{ex2}
   The following holds. \\
   (1)
   There are Hausdorff H-closed spaces of arbitrarily high Lindel\"of number. \\
   (2) There is a Hausdorff non-compact first countable space satisfying $\hu{\OO}{\Delta}$.\\
   (3)
   There is a first countable Lindel\"of H-closed Hausdorff space which does not satisfy $\hu{\OO}{\OO}$.
\end{prop}
\proof
The three examples are of the form $\wh{X}(\tau,\rho)$; 
Hausdorffness, H-closedness, Lindel\"ofness and first countability in (2) and (3) all follow from Lemma \ref{lemmawhX}.
\\
(1) This is well known, but let us give an example anyway.
    Take $X$ to be the ordinal $\kappa+1$, $\tau$ the order topology (which makes it compact) and $\rho$ the discrete topology.
    Then $L(\wh{X}(\tau,\rho))=\kappa$.\\
(2) Take $\kappa=\omega$ in (1). Then for each $\alpha\in\omega+1$, 
    $\{\alpha\}\times[0,1]$ is homeomorphic to $[0,1]$, so $\wh{X}(\tau,\rho)$ is a $\sigma$-compact space and thus satisfies $\hu{\OO}{\OO}$.
    We apply Lemma \ref{cbhu} to obtain $\hu{\OO}{\Delta}$. Of course, $\omega+1$ is first countable in the order topology.\\
(3) Take $X$ to be $[0,1]$, $\tau$ its usual topology, while $\rho$ is the coarsest refining of $\tau$ that makes $\Q\cap[0,1]$ clopen and discrete.
Thus, a $\rho$-open set is the union of (i) some subset of $\Q$ and (ii) $U-\Q$ with $U$ open for the usual topology.
Denote as usual the irrational numbers by $\mathbb{P}$.
It is well known that $\mathbb{P}\cap[0,1]$
is homeomorphic to the product space $\omega^\omega$ and does not satisfy $\hu{\OO}{\OO}$. Indeed,
fix a homeomorphism  $h:\omega^\omega\to\mathbb{P}\cap[0,1]$.
One way to easily obtain a sequence $\langle\mathcal{U}_n\rangle$ of covers of $\mathbb{P}\cap[0,1]$ violating $\hu{\OO}{\OO}$
is to set $\mathcal{U}_n=\{h(\pi_n^{-1}(\{m\}))\,:\,m\in\omega\}$
where $\pi_n$ is the projection on the $n$-th coordinate. 
Then the sequence of covers $\mathcal{W}_n=\{(U \cup (\Q\cap[0,1]))\times\{0\} \sqcup [0,1]\times (0,1]\,:\,U\in\mathcal{U}_n\}$
shows that $\wh{X}(\tau,\rho)$ does not satisfy $\hu{\OO}{\OO}$.
\endproof


\subsection*{The Property $\hu{\Delta}{\OO}$}\label{sec3}

We denote by
$\Delta_1$ the collection of open covers with at least one dense member.
First, some easy facts.

\begin{lemma} \label{lemmaequiv}
  Let $X$ be a space. The items below are equivalent:\\
  (a) $X$ satisfies $\hu{\Delta}{\OO}$,\\
  (b) $X$ satisfies $\hu{\Delta}{\Delta}$,\\
  (c) $X$ satisfies $\hu{\Delta_1}{\OO}$,\\
  (d) any closed subset of $X$ satisfies $\hu{\Delta}{\OO}$,\\
  (e) any closed nowhere dense subset of $X$ satisfies $\hu{\OO}{\OO}$.
\end{lemma}
\proof
  (c) $\rightarrow$ (a) $\leftrightarrow$ (b) are immediate, and (d) $\rightarrow$ (a) as well.\\
  (a) $\rightarrow$ (d) Let $Y\subset X$ be closed. Any od-cover of $Y$ yields an od-cover of $X$ by taking the union of the members with $X-Y$, and
     the result follows.\\
  (a) $\rightarrow$ (e) 
     Let $Y\subset X$ be closed and nowhere dense. If $Y$ does not satisfy $\hu{\OO}{\OO}$ take a sequence of covers 
     $\langle\mathcal{U}_n\rangle$ witnessing this fact.
     Set $\mathcal{V}_n =\{ U\cup (X-Y)\,:\,U\in \mathcal{U}_n\}$. Then $\langle\mathcal{V}_n\rangle$ witnesses that 
     $X$ does not satisfy $\hu{\Delta}{\OO}$.\\
  (e) $\rightarrow$ (c) 
     Let $\langle\mathcal{U}_n\rangle$ be a sequence of covers of $X$ ($n\in\omega$) such that some $U\in\mathcal{U}_0$ is dense in $X$.
     Set $\mathcal{F}_0 = \{U\}$.
     Since $X-U$ is closed and nowhere dense, there are finite $\mathcal{F}_n\subset\mathcal{U}_n$, $n\ge 1$,
     such that $\bigcup_{n\ge 1}\cup\mathcal{F}_n\supset X-U$. Then $\bigcup_{n\ge 0}\cup\mathcal{F}_n = X$.     
\endproof

The following proposition settles most of the classical cases (such as sets of reals).

\begin{prop}\label{propdeltao}
  Let $X$ be a separable space. Then $X$ satisfies $\hu{\OO}{\OO}$ if and only if $X$ satisfies $\hu{\Delta}{\OO}$.
\end{prop}
\proof
  One direction is trivial, so let us assume that $X$ satisfies $\hu{\Delta}{\OO}$.
  Let $D=\{d_i\,:\,i\in\omega\}$ be dense in $X$. 
  Given a sequence of open covers $\langle\mathcal{U}_i\,:\,i\in\omega\rangle$,
  take $V_{2i}\in\mathcal{U}_{2i}$ containing $d_i$ and set $\mathcal{F}_{2i}=\{V_{2i}\}$. 
  Since $V=\cup_{i\in\omega}V_{2i}$ contains $D$, 
  $X-V$ is closed and nowhere dense and satisfies $\hu{\OO}{\OO}$ by Lemma \ref{lemmaequiv}.
  Hence
  there are finite $\mathcal{F}_{2i+1}\subset\mathcal{U}_{2i+1}$ such that $\bigcup_{i\in\omega}\cup\mathcal{F}_{2i+1}\supset X-V$.
  Then $\bigcup_{i\in\omega}\cup\mathcal{F}_i = X$.
\endproof

Of course, od-compact spaces trivially satisfy $\hu{\Delta}{\OO}$. Any non-Lindel\"of such space (for instance: an uncountable discrete space) is a trivial
example of a space satisfying $\hu{\Delta}{\OO}$ but not $\hu{\OO}{\OO}$.
But we do not know the answer to the following question:

\begin{q}
   Is there a Lindel\"of non-od-compact space satisfying $\hu{\Delta}{\OO}$ but not $\hu{\OO}{\OO}$?
\end{q}

Another question, inspired by Theorem \ref{ehouais}:

\begin{q}
   Let $X$ be a space and $D\subset X$ the subspace of its isolated points. 
   Does the following equivalence hold:
   $X$ satisfies $\hu{\Delta}{\OO}$ $\longleftrightarrow$ $X-D$ satisfies $\hu{\OO}{\OO}$~?
\end{q}

{\footnotesize
\vskip .4cm
\noindent
Mathieu Baillif \\
Haute \'ecole du paysage, d'ing\'enierie et d'architecture (HEPIA) \\
Ing\'enierie des technologies de l'information TIC\\
Gen\`eve -- Suisse
}

\end{document}